\documentclass[a4paper,11pt,reqno]{amsart}

\usepackage[utf8]{inputenc}
\usepackage[british]{babel}
\usepackage[margin=2.5cm]{geometry}
\usepackage[shortlabels]{enumitem}
\usepackage{amsmath, amsthm, amssymb, amsfonts}
\usepackage{mathtools}
\usepackage{thmtools, thm-restate}
\usepackage{hyperref}
\usepackage{stmaryrd}
\usepackage{bbm}
\usepackage{tikz}
\usepackage{verbatim}

\newcommand{\cA}{\ensuremath{\mathcal A}}

\newcommand{\cD}{\ensuremath{\mathcal D}}
\newcommand{\cE}{\ensuremath{\mathcal E}}
\newcommand{\cF}{\ensuremath{\mathcal F}}

\declaretheorem[parent=section, name=Theorem]{thm}
\declaretheorem[sibling=thm]{lemma}

\declaretheorem[sibling=thm]{proposition}
\declaretheorem[sibling=thm]{claim}
\declaretheorem[sibling=thm]{conjecture}

\definecolor{RoyalAzure}{rgb}{0.0, 0.22, 0.66}
\definecolor{ForestGreen}{rgb}{0.13, 0.55, 0.13}

\renewcommand{\leq}{\leqslant}
\renewcommand{\geq}{\geqslant}

\newcommand*{\NN}{\mathbb{N}}

\DeclareMathOperator*{\Var}{Var}
\DeclareMathOperator*{\Cov}{Cov}
\DeclareMathOperator{\e}{\mathbb{E}}
\DeclareMathOperator*{\E}{\mathbb{E}}
\renewcommand{\Pr}{\mathbb{P}}

\newcommand{\Vh}{\ensuremath{V_\mathrm{heavy}}}
\newcommand{\Vl}{\ensuremath{V_\mathrm{light}}}
\newcommand{\dmed}{\ensuremath{d_\mathrm{med}}}
\DeclareMathOperator{\ind}{ind}
\DeclareMathOperator{\geom}{Geom}

\newcommand\subsetsim{\mathrel{\substack{
  \textstyle\subset\\[-0.2ex]\textstyle\sim}}}

\hypersetup{
  colorlinks,
  linkcolor={red!60!black},
  citecolor={green!50!black},
  urlcolor={blue!80!black}
}

\title{The edge-statistics conjecture for $\ell\ll k^{6/5}$}

\author{Anders Martinsson}
\address{Anders Martinsson, Department of Computer Science, ETH Z\"{u}rich, 8092 Z\"{u}rich, Switzerland}
\email{anders.martinsson@inf.ethz.ch}

\author{Frank Mousset}
\address{Frank Mousset, School of Mathematical Sciences, Tel Aviv University, Tel Aviv 6997801, Israel}
\email{moussetfrank@gmail.com}

\author{Andreas Noever}
\address{Andreas Noever, Department of Computer Science, ETH Z\"{u}rich, 8092 Z\"{u}rich, Switzerland}
\email{andreas.noever@gmail.com}

\author{Milo\v{s} Truji\'{c}}
\address{Milo\v{s} Truj\'{c}, Department of Computer Science, ETH Z\"{u}rich, 8092 Z\"{u}rich, Switzerland}
\email{mtrujic@inf.ethz.ch}

\date{\today}

\thanks{Research supported by ISF grants 1028/16 and 1147/14, and ERC Starting Grant 633509 (FM), and by grant no.\ 200021 169242 of the Swiss National Science Foundation (MT)}
\thanks{Part of this work has been completed at a workshop of the research group of Angelika Steger in Buchboden in July 2018}

\begin{document}

\begin{abstract}
  Let $k$ and $\ell$ be positive integers. We prove that if $1 \leq \ell \leq o_k(k^{6/5})$, then in every large enough graph $G$, the fraction of $k$-vertex subsets that induce exactly $\ell$ edges is at most $1/e + o_k(1)$. Together with a recent result of Kwan, Sudakov, and Tran, this settles a conjecture of Alon, Hefetz, Krivelevich, and Tyomkyn.
\end{abstract}

\maketitle

\section{Introduction}

Given a graph $G$ and some $k\in \NN$, let us write $X_{G,k}$ for the
random variable corresponding to the number of edges induced by a subset
$A\subseteq V(G)$
chosen uniformly at random among all subsets of size $k$.
Define $I(n,k,\ell) := \max{\{ \Pr[X_{G,k} =
\ell] \colon v(G) = n\}}$, the maximum probability of $X_{G,k}=\ell$ among all
$n$-vertex graphs $G$. A standard averaging argument shows that the function
$I(n,k,\ell)$ is decreasing in $n$, which implies that the limit \[
  \ind(k,\ell) := \lim_{n\to\infty} I(n,k,\ell)\] exists. Observe that
$\ind(k, \ell) = \ind(k, \binom{k}{2} - \ell)$. By
considering the empty/complete graphs on $n$ vertices, it is moreover easy to see that
$\ind(k,0) = \ind (k,\binom{k}{2}) = 1$, for all $k$. However, once
we exclude the cases $\ell \in
\{0, \binom{k}{2}\}$, it is sensible to suspect that $\ind(k, \ell)$ becomes much
smaller. For example, the quantitative version of Ramsey's theorem
implies that if $G$ is sufficiently large, then there is a positive
probability that $A$ is either a clique or an independent set, which shows
that $\ind(k,\ell)< 1$ for all $\ell \notin \{0,\binom{k}{2}\}$.

The function $\ind(k, \ell)$ was introduced by Alon, Hefetz, Krivelevich, and
Tyomkyn~\cite{alon2018edge}, motivated by a connection to the notion of graph
inducibility introduced earlier by Pippinger and
Golumbic~\cite{pippenger1975inducibility} (and which has recently become a
rather popular topic, see for example \cite{balogh2016maximum,
hatami2014inducibility, hefetz2018inducibility, kral2019bound,
yuster2019exact}). In \cite{alon2018edge}, Alon, Hefetz, Krivelevich, and
Tyomkyn advanced three conjectures concerning the asymptotics of the function
$\ind(k, \ell)$ as $k \to \infty$.
\begin{conjecture}\label{conj1}
  For all $k,\ell\in \NN$ with $0<\ell<\binom{k}{2}$, we have $\ind(k,\ell) \leq
  1/e + o_k(1)$.
\end{conjecture}

\begin{conjecture}\label{conj2}
  For all $k,\ell\in \NN$ with $\min\{\ell,\binom{k}{2}-\ell\} = \omega_k(k)$,
  we have $\ind(k,\ell)=o_k(1)$.
\end{conjecture}

\begin{conjecture}\label{conj3}
  For all $k,\ell\in \NN$ with $\min\{\ell,\binom{k}{2}-\ell\} = \Omega_k(k^2)$,
  we have $\ind(k,\ell)=O_k(k^{-1/2})$.
\end{conjecture}

\noindent
Here, the subscript $k$ indicates that the asymptotic notation is understood as
$k\to \infty$; for example, $o_k(1)$ denotes a function of $k$ tending to zero as
$k\to \infty$.
Several partial results on all three conjectures are given
in~\cite{alon2018edge}.

Note that Conjecture~\ref{conj2} implies Conjecture~\ref{conj1} in the range
where $\ell = \omega_k(k)$. Very recently, Kwan, Sudakov, and
Tran~\cite{kwan2019anticoncentration} gave a proof of Conjecture~\ref{conj2}
and showed that Conjecture~\ref{conj3} holds up to a polylogarithmic factor in
$k$. The purpose of this paper is to give a proof of Conjecture~\ref{conj1} for
all $1 \leq \ell \leq o_k(k^{6/5})$. Together with the result
of~\cite{kwan2019anticoncentration}, this result thus implies
Conjecture~\ref{conj1} for all $\ell$.

\begin{thm}\label{thm:main-theorem}
  For every $\ell= \ell(k)\in \NN$ such that $1\leq \ell \leq o_k(k^{6/5})$, we
  have \[ \ind(k,\ell) \leq 1/e + o_k(1). \]
\end{thm}

Even more recently, and independently of our own work, Fox and
Sauermann~\cite{fox2018comp} also gave a proof of Conjecture~\ref{conj1}.
The proof given here has the advantage that it is considerably shorter
than the one in~\cite{fox2018comp}. However, \cite{fox2018comp} contains
some stronger bounds in certain ranges of $\ell$ (e.g., it is shown
that in fact $\ind(k,\ell) = o_k(1)$ when $\omega_k(1) \leq \ell\leq o_k(k)$), as well as results
for the analogous problem in hypergraphs and other related results.

As noted in~\cite{alon2018edge}, the upper bound $1/e+o_k(1)$ in
Theorem~\ref{thm:main-theorem} is tight for example when $\ell=1$,
as can be seen by considering a random graph $G_{n,p}$ where $p =
1/\binom{k}{2}$. Similarly, the upper bound is tight for
$\ell=k-1$, as evidenced by
the complete
bipartite graph with parts of size $n/k$ and $(k-1)n/k$.
It would be interesting to know whether
the bound given by Theorem~\ref{thm:main-theorem} is tight for some values of
$\ell$ besides $1$ and $k-1$.

\section{A short proof for the case \texorpdfstring{$\ell = o_k(k)$}{l = o(k)}}\label{sec:short}
Before presenting the full proof of Theorem \ref{thm:main-theorem}, we give a short and self-contained proof for the case when $\ell=o_k(k)$.
\begin{proposition}\label{prop:small-ell}
  For every $\ell= \ell(k)\in \NN$ such that $1\leq \ell
  \leq o_k(k)$, we have \[ \ind(k,\ell) \leq 1/e + o_k(1). \]
\end{proposition}

\begin{proof}
  Choose $k$ and $\ell$ as in the statement and assume
  that $n=n(k)$ is sufficiently large.
  Let $G = (V,E)$ be a graph on $n$ vertices
  and let $\mathbf{v} = (v_1,v_2,\dotsc)$ be an infinite sequence of
  vertices chosen uniformly at random from $V^\NN$. We inductively colour the
  vertices in $\mathbf{v}$ with two colours, {\em black} and {\em green}, as follows:
  \begin{enumerate}[(1)]
    \item $v_1$ is black;
    \item $v_i$ is green if and only if the graph induced by $v_i$ and the
      \emph{black} vertices $v_j$ with $j<i$ contains at least $\ell$ edges;
      otherwise, $v_i$ is black.
  \end{enumerate}

  Set $L = L(\mathbf v) := \min{\{i\geq 1\colon \text{there are $k-1$ black
  vertices among $v_1,\dotsc,v_i$}\}}$ and $L := \infty$ if there are fewer
  than $k-1$ black vertices in $\mathbf v$. We then define $Y_{G,k} = Y_{G,k}(\mathbf v)$ as
  the random variable corresponding to the number of green vertices in the set
  $\{v_i\colon 1\leq i < L\}$.

  We first show that
  \begin{equation}\label{eq:couple}
    \Pr[X_{G,k} = \ell] \leq \Pr[Y_{G,k} = 1]+o_k(1).
  \end{equation}
  This can be seen as follows. Let $\tilde X_{G,k} =
  e(\{v_1,\dotsc,v_k\})$ and let $\cA$ be the event that
  $v_1,\dotsc,v_k$ are all distinct. If $n$ is sufficiently
  large given $k$ (i.e., $n = \omega(k^2)$), then
  $\Pr[\cA]=1-o_k(1)$. Thus
  \begin{equation}\label{eq:all-distinct}
    \Pr[X_{G,k} = \ell] = \Pr[\tilde X_{G,k} = \ell \mid \cA]
    \leq \Pr[\tilde X_{G,k} =
    \ell]/\Pr[\cA] \leq \Pr[\tilde X_{G,k} = \ell] + o_k(1).
  \end{equation}
  Next, since $\ell$ edges can span at most $2\ell$ vertices, it follows by
  symmetry that
  \begin{equation}\label{eq:last-not-bad}
    \Pr[\tilde X_{G,k} = \tilde X_{G,k-1} =\ell] \geq \Pr[\tilde
    X_{G,k} = \ell]\cdot \frac{k-2\ell}{k} \geq \Pr[\tilde
    X_{G,k} = \ell] - o_k(1),
  \end{equation}
  where the last inequality uses $\ell = o(k)$ (in fact, this is the only place where we use this assumption).
  Finally, and crucially, observe that $\tilde X_{G,k} = \tilde X_{G,{k-1}} =
  \ell$ implies $Y_{G,k} =1$:
  if $\tilde X_{G,k} = \tilde X_{G,k-1}\geq \ell$, then at least one green
  vertex must appear before the $(k-1)$-st black vertex, and if
  there is more than one
  such green vertex, then
  $\tilde X_{G,k}>\ell$.
  From this, together with \eqref{eq:all-distinct}
  and \eqref{eq:last-not-bad}, it follows that
  $\Pr[X_{G,k} = \ell] \leq \Pr[Y_{G,k} = 1] + o_k(1)$,
  as claimed. Therefore, it suffices to show that
  $\Pr[Y_{G,k}=1]\leq 1/e$.

  Let $\mathbf u = (u_1,\dotsc,u_{k-1})$ be a sequence of $k-1$ (not
  necessarily distinct) vertices of $G$. Let $U(\mathbf u)$ be the event that
  $u_1,\dotsc,u_{k-1}$ are the first $k-1$ black vertices in $\mathbf v$.
  Now observe that if $\Pr[U(\mathbf u)]$ is nonzero, then the conditional
  distribution of $Y_{G,k}$ given $U(\mathbf u)$ is given by the sum
  \[ \geom(p_1) + \geom(p_2) + \dotsb + \geom(p_{k-2}) \]
  of independent geometric distributions with parameters
  \[ p_i := \frac{1}{n} \big| \{v \in V \colon
  e(\{u_1,\dotsc,u_i,v\}) \geq \ell\} \big|. \]
  Indeed, suppose that we have chosen the vertices
  $v_1,v_2,\dotsc,v_t = u_i$ up to $u_i$.
  From then on, each vertex that we choose from the set $\{v \in V \colon
  e(\{u_1,\dotsc,u_i,v\}) \geq \ell\}$ is green, while the first vertex that we
  choose outside of this set is the next black vertex $u_{i+1}$. It follows that
  \[
    \Pr[Y_{G,k}=1\mid U(\mathbf u)]
    = \sum_{i=1}^{k-2} p_i \prod_{j=1}^{k-2} (1-p_j)
    \leq \sum_{i=1}^{k-2} p_i \cdot e^{-\sum_{j=1}^{k-2} p_j}
    \leq 1/e,
    \]
  using that $f(x) = xe^{-x}$ is maximised for $x=1$. Since this is true for
  every relevant choice of $\mathbf u$, we also have
  $\Pr[Y_{G, k} = 1]\leq 1/e$ unconditionally. The proposition then follows using \eqref{eq:couple}.
\end{proof}

\section{Proof of Theorem~\ref{thm:main-theorem}}

We use the following simple facts about hypergeometric random variables.
\begin{lemma}\label{lem:hypergeom}
  Let $X$ be hypergeometric random variable counting the number of successes
  obtained when sampling $m$ elements without replacement from a population of
  size $N$ containing $Np$ successes.
  Assume $m^2/N \to 0$ and $m\to\infty$.
  If $mp\to \lambda<\infty$, then
  \[ \max_{i} \left|\Pr[X=i] - \frac{\lambda^ie^{-\lambda}}{i!}\right| \to
  0,\]
  where the maximum is taken over all nonnegative integers.
  On the other hand, if $mp(1-p) \to \infty$, then
  $\max_i\Pr[X=i]\to 0$.
\end{lemma}

\begin{proof}
  Let $Y$ be defined in the same way as $X$ except that the $m$ elements are
  sampled with replacement. Let $\cA$ be the event that all $m$ elements
  are distinct in this experiment.
  Then we have $\Pr[X=i] = \Pr[Y=i \mid \cA]$.
  The assumption $m^2/N\to 0$ implies that $\Pr[\cA]\to 1$ and hence
  $\max_i \big| \Pr[X=i] - \Pr[Y= i] \big| \to 0$. Note that $\Pr[Y=i]
  = p^i (1-p)^{m-i}\binom{m}{i}$. Then the first assertion follows
  from the usual Poisson approximation to the binomial distribution.
  Similarly, the second assertion follows from the
  de Moivre--Laplace theorem.
\end{proof}

Let $k$ and $\ell$ be such that $1\leq \ell \leq o_k(k^{6/5})$
and
assume that $G$ is graph with $n$ vertices, where we assume
that $n=n(k)$ is
sufficiently large to support
our arguments.
We always interpret asymptotic statements
as $k\to \infty$, and thus omit the subscript $k$ in the asymptotic
notation from now on. We say that an event holds {\em with high probability} (w.h.p.\ for short) if the probability that it holds approaches $1$ as $k\to\infty$.

For two events $\cE=\cE(k)$ and $\cF=\cF(k)$ (which can thus also
depend on $\ell$, $G$, and $n$), we say that \emph{$\cE$ is
essentially contained in
$\cF$}, and write $\cE \subsetsim \cF$, if $\Pr(\cE \setminus \cF) =o(1)$.

As in the introduction, let $A$ denote a uniformly random subset of $V(G)$ of size $k$
and set $X_{G,k} = e(A)$. Throughout the proof, we let $\cE$
denote the event that $X_{G,k} = \ell$.

Observe that it is enough to show that $\cE$ is essentially contained in an event
of probability $1/e+o(1)$. To define this event, let first
$(w_k)_{k\geq 1}$ be a sequence of positive real numbers that goes to infinity
at a sufficiently slow rate and, for every integer $d\geq 0$, define the event
\[ \cD_d :=
  \{\text{all but at most $w_k\sqrt{\ell}$ vertices in $A$ have degree $d$ in
$G[A]$} \}.
\]
In particular, we choose $w_k$ such that $w_k\sqrt{\ell} =o(k)$. Our main goal is to show that there exists some deterministic value $d=d(G,k,\ell)$ such that $\cE\subsetsim \cE\cap \cD_d$.
This is sufficient by the following claim.

\begin{claim}\label{cl:Dd}
  For every $d\geq 0$, we have
  \[ \Pr[\cE\cap \cD_d] \leq 1/e + o(1). \]
\end{claim}

\begin{proof}
  Assume first that $d\geq 1$.
  Let $v$ be a vertex chosen uniformly at random among the vertices
  in $A$. Since $\cD_d$ implies that all but $o(k)$ vertices
  of $A$ have degree $d$ in $G[A]$, we have $\Pr[e(v,A) = d \mid \cD_d] =
  1-o(1)$ and thus
  \[ \Pr[\cE\cap \cD_d] \leq \Pr[\cD_d] \leq (1+o(1))\Pr[e(v, A) = d\text{ and
  } \cD_d] \leq \Pr[e(v, A) = d] +o(1). \]
  Note also that we have $\Pr[\cE \cap
  \cD_d] = 0$ unless $kd \leq 3\ell$ and so
  using
  $\ell = o(k^{6/5})$ we can assume
  $d=o(k^{1/5})$.
  Note next that we can generate the pair $(v,A)$ by choosing first a
  uniformly random vertex $v$ in $V(G)$ and then choosing the remaining $k-1$
  vertices of $A$ uniformly among the $(k-1)$-element subsets of $V(G)\setminus
  \{v\}$. In particular, if we fix the choice of $v$, then $e(v,A)$ follows
  a hypergeometric distribution with sample size $k-1$ and a population of
  size $n-1$ comprising $d_G(v)$ successes.
  If $d_G(v) \geq n/2$ then it
  follows from $d=o(k)$ and Markov's inequality that $e(v,A)>d$ with
  probability $1-o(1)$.
  On the other hand, if $d_G(v)\leq n/2$, then Lemma~\ref{lem:hypergeom}
  implies that we
  either have
  $\Pr[e(v,A)=d]
  =o(1)$ or
  \[ \Pr[e(v,A) = d] =
  \frac{\lambda^d e^{-\lambda}}{d!} +o(1)\]
  for some $\lambda \geq 0$. Optimising the value of $\lambda$, we see that
  \[ \frac{\lambda^d e^{-\lambda}}{d!}
  \leq \frac{d^d e^{-d}}{d!}\leq 1/e, \]
  where the last inequality uses $d\geq 1$.

  Suppose next that $d=0$. In this case we can proceed as in the
  proof of Proposition~\ref{prop:small-ell}, if, instead of
  \eqref{eq:couple}, we can show that
  \begin{equation}\label{eq:couple2}
    \Pr[\cE\cap \cD_0] \leq \Pr[Y_{G,k} = 1]+o(1).
  \end{equation}
  Assume the process is the same as in
  Proposition~\ref{prop:small-ell} and that $Y_{G, k}$, $\tilde
  X_{G, k}$, and $\cA$ are defined in the same way. Then
  \eqref{eq:couple2} can be seen as follows.
  Let $\tilde \cD_0$ be the event that all but at most $w_k\sqrt{\ell} = o(k)$
  of the vertices $v_1,\dotsc, v_k$ are isolated in $G[\{v_1,\dotsc,v_k\}]$.
   We have
  \[
    \Pr[\cE\cap \cD_0] = \Pr[\tilde \cD_0 \text{ and } {\tilde X_{G,k} = \ell} \mid \cA] \leq \Pr[\tilde \cD_0 \text{ and } {\tilde X_{G,k} = \ell}] / \Pr[\cA] \leq \Pr[\tilde \cD_0 \text{ and } {\tilde X_{G,k} = \ell}] + o(1).
  \]
  Since each permutation of $v_1,\dotsc,v_k$ is equally likely, we further have
  \[
    \Pr[\tilde X_{G,k} = \tilde X_{G,k-1} =\ell] \geq \Pr[\tilde \cD_0 \text{ and } {\tilde X_{G,k} = \ell}]-\frac{w_k\sqrt{\ell}}{k},
  \]
  where the error term in the right hand side is $o(1)$ provided
  $w_k$ increases slowly enough. As $\tilde X_{G,k} =
  \tilde X_{G,k-1} = \ell$ implies $Y_{G,k}=1$ deterministically, the
  proof of \eqref{eq:couple2} is complete.
\end{proof}

It remains to show that there is some $d = d(G, k, \ell)$ such that
$\cE\subsetsim \cD_d$. We do this over a series of claims. First,
let us define the event
\[ \cD_* = \bigcup_{d\geq 0} \cD_d =
\{ \text{all but at most $w_k \sqrt\ell$ vertices in $A$ have the same degree
in $G[A]$}\}. \]
The first claim we need is the following:

\begin{claim}\label{claim:dstar}
  We have $\cE \subsetsim \cD_*$.
\end{claim}

The somewhat technical proof of Claim~\ref{claim:dstar} is
deferred to the end of the paper. With this claim at hand, we
continue with the proof of the theorem. We partition the
vertices of $G$ into two sets:
\begin{itemize}
  \item the \emph{heavy} vertices $\Vh := \{v\in V(G)
    \colon \deg_G(v) \geq n\ell^{1/3}/k\}$;
  \item the \emph{light} vertices $\Vl :=
    \{v\in V(G) \colon \deg_G(v) < n\ell^{1/3}/k\}$.
\end{itemize}
We first show that we can assume that there are not too
many heavy vertices.

\begin{claim}\label{claim:heavy}
  Assume that $\ell =\omega(1)$ and that $G$ contains more than $5\ell^{2/3}
  n/k$ heavy vertices. Then
  \[ \Pr[\cE] = o(1).  \]
\end{claim}
\begin{proof}
  We generate $A$ by first choosing a random set $A_1$ of size
  $k/2$ and then choosing another random set $A_2\subseteq
  V(G)\setminus A_1$ of size
  $k/2$. We have
  \[ \E \big[ |A_1 \cap \Vh| \big]
  = \frac{ |\Vh|k}{2n} \geq
  \frac{5\ell^{2/3}}{2} = \omega(1). \]
  In particular, the Chernoff bounds for the hypergeometric
  distribution imply that w.h.p.\ $A_1$ contains at least
  $2.49\ell^{2/3}$ heavy vertices. Expose the set $A_1$ and assume that this is the case. Then every (fixed) vertex $v \in A_1 \cap \Vh$ satisfies
  $\e[e(v, A_2)] \geq (1-o(1)) \ell^{1/3}/2$. Hence, again by
  the Chernoff bounds and a union bound over an arbitrary set $S \subseteq A_1 \cap \Vh$ of size $2.49 \ell^{2/3}$, we get
  \[
    \Pr[\exists v \in S \colon e(v,A_2) < 0.49\ell^{1/3}] = o(1).
  \]
  In particular, the union $A = A_1\cup A_2$ w.h.p.\ contains at least
  \[ 2.49\ell^{2/3}\cdot 0.49\ell^{1/3}> \ell \]
  edges of $G$, implying $\Pr[\cE] = o(1)$.
\end{proof}

\begin{claim} \label{claim:hhll}
  Let $Z := \sum_{v \in A \cap \Vl} e(v, A)$.
  Assume that $\ell= \omega(1)$.
  Then either
  \[ \Pr[\cE] =o(1) \]
  or
  \[ \Var\left[X_{G, k} - Z \right] \leq 30\ell^{5/3}. \]
\end{claim}
\begin{proof}
  Let $H := e(A\cap \Vh)$ and $L := e(A\cap \Vl)$ and observe that $X_{G,k}-Z
  = H - L$. Using the elementary inequality $(a - b)^2 \leq 2 a^2 + 2b^2$, we have
  \[
    \Var[X_{G, k} - Z] = \Var[H - L] \leq 2\Var[H] + 2\Var[L].
  \]

  For any edge $e\in G$, let $X_e$ denote the indicator random variable for the event that both
  endpoints of $e$ are contained in $A$. We have
  \[
    \Var[H] = \sum_{e \in G[\Vh]} \sum_{f \in G[\Vh]} \Cov[X_e, X_f]
  \]
  and
  \[
    \Var[L] = \sum_{e \in G[\Vl]} \sum_{f \in G[\Vl]} \Cov[X_e, X_f].
  \]

  For each of these sums, an elementary calculation shows that $\Cov[X_e, X_f] \leq
  0$ if $e$ and $f$ do not have a common endpoint. On the other hand, if $e$
  and $f$ intersect in exactly one endpoint, we have
  $\Cov[X_e, X_f]\leq \e[X_eX_f]\leq \e[X_e] \cdot (k/n)$. Lastly, we have
  $\Cov[X_e, X_e] = \Var[X_e] \leq \e[X_e]$.

  Let $\mu_1:=\e[H]$ and $\mu_2:=\e[L]$.
  Since we may assume $|\Vh|\leq 5\ell^{2/3} n/k$ (as otherwise
  Claim~\ref{claim:heavy} implies $\Pr[\cE]=o(1)$), we then obtain
  \[\Var[H] \leq e(\Vh) \cdot \e[X_e] + 2 \cdot e(\Vh) \cdot \frac{5\ell^{2/3} n}{k}
  \cdot \e[X_e] \cdot \frac{k}{n}
  \leq
  (1+o(1)) \mu_1 \cdot 10\ell^{2/3}.
  \]
  Similarly, using the fact that every light vertex has degree at most $n\ell^{1/3}/k$, we get
  \[ \Var[L] \leq
  e(\Vl) \cdot \e[X_e] + 2 \cdot e(\Vl) \cdot \frac{n\ell^{1/3}}{k} \cdot \e[X_e]\cdot \frac{k}{n}
  \leq
  (1+o(1))\mu_2 \cdot 2\ell^{1/3}.
  \]

  If either of $\mu_1$ or $\mu_2$ is greater than, say, $1.01\ell$, then
  by Chebyshev's inequality, the corresponding random variable $H$ or $L$
  is concentrated around its expectation, which (since $H, L\leq X_{G,k}$)
  would imply that $\Pr[X_{G,k} = \ell] = o(1)$.
  Otherwise, if $\mu_1,\mu_2\leq 1.01\ell$, we obtain the
  desired upper bound on $\Var[X_{G, k} - Z]$.
\end{proof}

\begin{claim}\label{claim:endgame}
  Assume that $\ell = \omega(\log^{3} k)$. Then 
  there exists some deterministic $d = d(G,k,\ell)$ such that
  $\cE\subsetsim \cD_d$.
\end{claim}
\begin{proof}
  By
  Claim~\ref{claim:heavy}, we can assume that there are at most
  $5\ell^{2/3}n/k$ heavy vertices in $G$, since otherwise $\Pr[\cE] = o(1)$
  and then $\cE\subsetsim \cD_0$ (say) holds trivially.

  As in the statement of Claim~\ref{claim:hhll}, let $Z := \sum_{v\in A\cap
  \Vl} e(v, A)$. Again, since we are done when $\Pr[\cE] = o(1)$,
  we can assume that
  \begin{equation}
    \label{eq:var}
    \Var[X_{G, k} - Z] \leq 30\ell^{5/3},
  \end{equation}
  using Claim~\ref{claim:hhll}.

  We denote by $D$ the random variable corresponding to the most frequent
  degree in $G[A]$ (with ties broken arbitrarily). We first show that $\cE$ is
  essentially contained in each of the following events:
  \begin{itemize}
    \item $\cF_1:=\{$every $v\in A\cap \Vl$ satisfies $e(v, A)\leq 2 \ell^{1/3}\}$,
    \item $\cF_2:=\{$every $v\in A\cap \Vh$ satisfies $e(v, A) \geq
      \ell^{1/3}/2\}$,
    \item $\cF_3:=\{X_{G,k}=Z+\mu \pm w_k \ell^{5/6}\}$,
      where $\mu = \e[X_{G,k}-Z]$,
    \item $\cF_4:=\{Z = kD \pm 3w_k \ell^{5/6}\}$.
  \end{itemize}

  Since $\ell^{1/3} = \omega(\log k)$, the Chernoff bounds easily imply
  $\Pr[\cF_1\cap \cF_2] = 1-o(1)$, so
  $\cE\subsetsim \cF_1$ and $\cE\subsetsim \cF_2$ hold trivially. For $\cF_3$,
  note that using \eqref{eq:var},
  Chebyshev's inequality gives $\Pr[\overline{\cF_3}] \leq O(1/w_k^2)
  =o(1)$, thus we have $\cE\subsetsim \cF_3$ as well.

  By Claim~\ref{claim:dstar}, we further know that $\cE\subsetsim \cD_*$, and therefore $\cE \subsetsim \cE \cap \cD_* \cap \cF_1 \cap \cF_2$.
  To prove that $\cE\subsetsim \cF_4$, it is thus enough to show
  that $\cE \cap \cD_* \cap \cF_1\cap \cF_2 \subseteq \cF_4$
  (note that this is a deterministic statement).
  So assume that $\cE \cap \cD_* \cap \cF_1\cap \cF_2$ holds.
  Since $D$ is the most frequent degree in $G[A]$, we see that
  $\cE\cap \cD_*$ implies $\ell = X_{G,k}\geq (k-o(k))D/2
  \geq kD/3$ for
  all sufficiently large $k$. As $\cF_2$ implies that every
  heavy vertex $v\in A$ satisfies $e(v, A)\geq \ell^{1/3}/2 \gg \ell/k$ (recall, $\ell = o(k^{6/5})$),
  all of the at least $k-w_k\sqrt{\ell}$
  vertices $v\in A$ with
  $e(v, A)=D$ are light. It follows that
  \[ (k-w_k\sqrt\ell) D \leq Z \leq kD + 2w_k
  \ell^{5/6}, \]
  where the upper bound is implied by $\cF_1$.
  Therefore, using $D\leq 3\ell/k$,
  \[ kD - w_k\sqrt{\ell} \cdot 3 \ell/k \leq Z \leq kD
  + 2w_k \ell^{5/6}. \]
  Since $\ell=o(k^{3/2})$, we have
  $\ell^{3/2} /k =o(\ell^{5/6})$, so
  the above implies
  $\cF_4$. It follows that $\cE\subsetsim \cF_4$.

  Finally, note that $\cE \cap \cF_3 \cap \cF_4$ gives
  \[ D = \frac{\ell-\mu}{k} \pm \frac{w_k}{k}\cdot O(
  \ell^{5/6}). \]
  By letting $w_k$ be a sufficiently slowly diverging function, the error term
  in the right hand side is $o(1)$ (using in addition $\ell = o(k^{6/5})$), meaning there is only
  (at most)
  one possible integer value of $D$ that can satisfy this. Let $d$ be this value.  Then $\cE \subsetsim \cE
  \cap \cD_* \cap \cF_3 \cap \cF_4 \subseteq \cD_d$, as desired.
\end{proof}

Claims~\ref{cl:Dd} and \ref{claim:endgame} imply that we have $\Pr[X_{G,k}=\ell] \leq
1/e+o(1)$ for all $\omega(\log^3 k)\leq \ell \leq o(k^{6/5})$ (and we already proved
the case $1\leq \ell = o(k)$ in Section~\ref{sec:short}). Thus it
only remains to prove Claim~\ref{claim:dstar}.

\subsection{Proof of Claim~\ref{claim:dstar}}

We now give the missing proof of Claim~\ref{claim:dstar}. Let $m =
k/(w_k^{1/3}\sqrt{\ell})$. If $w_k$ diverges sufficiently slowly, and using
$\ell = o(k^{6/5})$, we have (say) $m \geq w_k$. Observe that we can generate
$A$ by first choosing a random set $S$ of size $k-m$ and then choosing a
random set $Q$ of size $m$ from the complement of $S$. In terms of this
process, we define the
following events:
\begin{itemize}
  \item $\cE_1 := \{ e(Q) = 0 \}$,
  \item $\cE_2 := \{ e(S) + \sum_{v\in Q} e(v, S) = \ell \}$,
  \item $\cE_3 := \{ \text{all but at most $w_k^{1/3}$ vertices in $Q$ have the
    same degree into $S$}\}$,
  \item $\cE_4 := \{ \text{all but at most $w_k^{1/3}$ vertices in $Q$ have the same degree in $A$}\}$.
\end{itemize}
We prove that $\cE$ is essentially contained in each of these events,
and then use this to conclude that $\cE\subsetsim \cD_*$.

We first prove that $\cE \subsetsim \cE_1$.
Since we can generate $Q$ by first generating $A$ and then choosing a random subset
$Q\subseteq A$ of size $m$, we have
\[
  \E[e(Q)\mid X_{G,k} = \ell] = \ell \cdot \binom{m}{2} / \binom{k}{2} =
  O(1/w_k^{2/3}),
\]
where the last inequality uses the definition of $m$. Therefore, by Markov's inequality,
\[
  \Pr[{X_{G, k} = \ell} \text{ and } {e(Q) \neq 0}] \leq \Pr(X_{G, k} = \ell) \cdot O(1/w_k^{2/3}) = o(1),
\]
so $\cE\subsetsim \cE_1$.

Having this, it follows directly from the definitions that $\cE
\subsetsim \cE \cap \cE_1 \subseteq \cE_2$.

Next, we show that $\cE_2 \subsetsim \cE_3$, which then implies $\cE\subsetsim
\cE_3$. Expose first only the set $S$ and let $\dmed$ be the median of $e(v,
S)$ over all $v \in V(G) \setminus S$. We consider two cases, depending on the
properties of the set $S$.

\textit{Case 1.} All but at most $w_k^{1/4} n/m$ vertices $v \in V(G) \setminus
S$ satisfy $e(v, S) = \dmed$.  Clearly, the expected number of vertices $v \in
Q$ for which $e(v, S) \neq \dmed$ is then at most $O(w_k^{1/4}) = o(w_k^{1/3})$.
Thus, by Markov's inequality, we have $\Pr[\cE_3] = 1-o(1)$, which implies
$\cE_2\subsetsim \cE_3$ in this case.

\textit{Case 2.} At least $w_k^{1/4} n/m$ vertices $v \in V(G) \setminus S$
satisfy $e(v, S) \neq \dmed$. We claim that in this case, we have $\Pr[\cE_2] =
o(1)$. We can assume that at least $w_k^{1/4} n/(2m)$ vertices $v \in V(G)
\setminus S$ satisfy, say, $e(v, S) > \dmed$ (the case in which at least
$w_k^{1/4} n/(2m)$ vertices $v \in V(G) \setminus S$ satisfy $e(v, S) < \dmed$
is analogous). Let us denote the number of such vertices by $t$
and let $N := |V(G)\setminus S| = n-k+m$.

Note that we can generate the set $Q$ in the following way. First, let
$v_1',v_2',\dotsc, v_{N-t}'$ be a random permutation of the vertices $v\in
V(G)\setminus S$ with $e(v, S) \leq \dmed$, and let
$v_1'',v_2'',\dotsc,v_t''$ be a random permutation of the vertices $v\in
V(G)\setminus S$ with $e(v, S) > \dmed$.  Let $I$ be the random
variable corresponding to the number of red balls one obtains when drawing $m$
balls without replacement from a population of size $N$ containing $N-t$
red balls and $t$ blue balls (in other words, let $I$ be a hypergeometric random
variable with these parameters). Finally, let \[ Q = \{ v_1', v_2', \ldots,
v_I', v_1'', v_2'', \ldots, v_{m - I}'' \}. \] Note that in this way, $Q$ is
really a uniformly random $m$-element subset of $V(G)\setminus S$.

Now, in order for $\cE_2$ to occur we need
\[
  \sum_{v \in Q} e(v, S) = \ell - e(S).
\]
Observe that for every fixed choice of the permutations
$v_1',v_2',\dotsc,
v_{N-t}'$
and
$v_1'',v_2'',\dotsc,v_t''$, there is at most one value of $I$
that achieves this.
However, since $I$ is a hypergeometric random variable with population size
$N$ and sample size $m$, which satisfy
$m^2/N\leq k^2/(n-k)=o(1)$ if $n = \omega(k^2)$,
and since
$t \geq w_k^{1/4} n/(2m) = \omega(N/m)$ and $t\leq
N/2$ (as $\dmed$ is a median) imply
$m(t/N)(1-t/N)=\omega(1)$,
it follows from Lemma~\ref{lem:hypergeom} that
$\Pr[I=i]=o(1)$.
Thus in this case, we have $\Pr[\cE_2]=o(1)$, from which
$\cE_2\subsetsim \cE_3$ follows trivially.

Having shown $\cE\subsetsim \cE_1$ and $\cE\subsetsim \cE_3$, it follows
easily from the definitions that $\cE \subsetsim \cE_1 \cap \cE_3 \subseteq
\cE_4$.

Lastly, we show that $\cE_4 \subsetsim \cD_*$, which completes the proof.
Suppose that $A$ is such that $\cD^*$ does not occur. We show that,
conditioning on this event (but leaving the subset $Q\subseteq A$ random), the
probability of $\cE_4$ is $o(1)$.
For this, let $d$ be the median degree in $G[A]$.
Then at least $w_k \sqrt\ell/2$ vertices have degree, say, larger than $d$ in $G[A]$ (the case where $w_k\sqrt\ell /2$ vertices
have degree smaller than $d$ is analogous).
Let $t$ be the number of such vertices in $A$ and let $X_t$ be the random
variable denoting the number of such vertices in $Q$ (which,
recall, is a random subset of $A$ of size $m$). Then
since $m= k/(w_k^{1/3} \sqrt\ell)$, we have
\[
  \E[X_t] = t \cdot \frac{m}{k} \geq w_k \cdot \frac{m \sqrt\ell}{2k} =
  w_k^{2/3}/2= \omega(1),
\]
and $\sigma(X_t) = O(\sqrt{t m/k})$. Therefore, by Chebyshev's inequality, w.h.p.\ we have $w_k^{1/3} \leq X_t$. On the other hand, as $t \leq k/2$
(recall, $d$ is a median), we
also have w.h.p.\ $X_t \leq (1/2 + o(1))m$. Since $w_k \ll m$,
these two inequalities imply that
there is no set of $m - w_k^{1/4}$ vertices in $Q$ which have the same degree in $A$. Consequently, $\Pr(\cE_4 \mid \overline\cD_*) = o(1)$, which implies
$\Pr(\cE_4 \setminus \cD_*) = o(1)$, as desired. \qed

\bibliography{refs}
\bibliographystyle{abbrv}

\end{document}